\documentclass[12pt]{amsart}
\usepackage{amscd,amssymb,verbatim}
\usepackage{amssymb,amsmath,amsthm,epsfig}

\newcommand{\N}{\mathbb{N}}

\theoremstyle{plain}
\newtheorem{Thm}{Theorem}[section]

\newtheorem{Lem}[Thm]{Lemma}
\newtheorem{Cor}[Thm]{Corollary}

\newtheorem*{qst}{Question}
\newtheorem{Def}[Thm]{Definition}

\theoremstyle{remark}

\begin{document}

\title{A property of strictly singular 1-1 operators}
\author{George~Androulakis, Per~Enflo}
\thanks{The research was partially supported by NSF}
%\date{December 24, 2001}
\maketitle

\noindent
{\bf Abstract} We prove that if $T$ is a  strictly singular 1-1
operator defined on an infinite dimensional Banach space $X$, then for
every infinite dimensional subspace $Y$ of $X$ there exists an
infinite dimensional subspace $Z$ of  $Y$ such that $Z$ contains
orbits of $T$ of every finite length and the restriction of $T$ on $Z$
is a compact operator.

\section{Introduction} \label{sec1}

An operator on an infinite dimensional Banach space is called {\em
  strictly singular} if it fails to be an isomorphism when it is
restricted to any infinite dimensional subspace (by ``operator'' we
will always mean a ``continuous linear map''). It is easy to see that
an operator $T$ on an infinite dimensional Banach space $X$ is
strictly singular if and only if for every infinite dimensional
subspace $Y$ of $X$ there exists an infinite dimensional subspace $Z$
of $Y$ such that the restriction of $T$ on $Z$, $T|_Z: Z \to X$, is a
compact operator. Moreover, $Z$ can be assumed to have a
basis. Compact operators are special examples of strictly singular
operators. If $1 \le p < q \le \infty$ then the inclusion map 
$i_{p,q}: \ell_p \to \ell_q$ is a strictly singular (non-compact)
operator. A {\em Hereditarily Indecomposable} (H.I.) Banach space is
an infinite dimensional space such that no subspace can be written as
a topological sum of two infinite dimensional subspaces. W.T. Gowers
and B. Maurey constructed the first example of an H.I. space
\cite{GM}. It is also proved in \cite{GM} that every operator on a
complex H.I. space can be written as a strictly singular perturbation
of a multiple of the identity. If $X$ is a complex H.I. space and $T$
is a strictly singular operator on $X$ then the spectrum of $T$
resembles the spectrum of a compact operator on a complex Banach
space: it is either the singleton $\{ 0 \}$  (i.e. $T$ is
quasi-nilpotent), or a sequence 
$\{ \lambda_n : n=1,2, \ldots \} \cup \{ 0 \}$ where $\lambda_n$ is an
eigenvalue of $T$ with finite multiplicity for all $n$, and
$(\lambda_n)_n$ converges to $0$, if it is an infinite sequence. It
was asked whether there exists an H.I. space $X$ which gives a
positive solution to the ``Identity plus Compact'' problem, namely,
every operator on $X$ is a compact perturbation of a multiple of the
identity. This question was answered in negative in \cite{AS} for the
H.I. space constructed in \cite{GM}, (for related results see
\cite{Ga}, \cite{G}, and \cite{AOST}). By \cite{ArS}, (or the more
general beautiful theorem of V. Lomonosov \cite{L}), if a Banach space
gives a positive solution to the ``Identity plus Compact'' problem, it
also gives a positive solution to the famous Invariant Subspace
Problem (I.S.P.). The I.S.P. asks whether there exists a separable
infinite dimensional Banach space on which every operator has a
non-trivial invariant subspace, (by ``non-trivial'' we mean
``different than $\{ 0 \}$ and the whole space''). It remains unknown
whether $\ell_2$ is a positive solution to the I.S.P.. Several
negative solutions to the I.S.P. are known \cite{En1}, \cite{En2}, \cite{R1},
\cite{R2}, \cite{R3}. In particular, there exists a strictly singular
operator with no non-trivial invariant subspace \cite{R4}. It is
unknown whether every strictly singular operator on a super-reflexive
Banach space has a non-trivial invariant subspace. Our main result
(Theorem \ref{main}) states that if $T$ is a strictly singular 1-1
operator on an infinite dimensional Banach space $X$, then for every
infinite dimensional Banach space $Y$ of $X$ there exists an infinite
dimensional Banach space $Z$ of $Y$ such that the restriction of $T$
on $Z$, $T|_Z:Z \to X$, is compact, and $Z$ contains orbits of $T$ of
every finite length (i.e. for every $n \in \N$ there exists 
$z_n \in Z$ such that 
$\{ z_n, Tz_n, T^2z_n, \ldots, T^nz_n \} \subset Z$). We raise the
following

\begin{qst}
Let $T$ be a quasi-nilpotent operator on a super-reflexive Banach
space  $X$, such
that for every infinite dimensional subspace $Y$ of $X$ there exists
an infinite dimensional subspace $Z$ of $Y$ such that $T|_Z:Z \to X$ is
compact and $Z$ contains orbits of $T$ of every finite length. Does
$T$ have a non-trivial invariant subspace? 
\end{qst}

By our main result, an affirmative answer to the above question would
give that every strictly singular, 1-1, quasi-nilpotent operator on a
super-reflexive Banach space has a non-trivial invariant subspace; in
particular, we would obtain that every operator on the super-reflexive
H.I. space constructed by V. Ferenczi \cite{F} has a non-trivial
invariant subspace, and thus the I.S.P. would be answered in
affirmative. 

\section{The main result} \label{sec2}

Our main result is

\begin{Thm} \label{main}
Let $T$ be a strictly singular 1-1 operator on an infinite
dimensional Banach space $X$. Then, for every infinite dimensional
subspace $Y$ of $X$ there exists an infinite dimensional subspace $Z$
of $Y$, such that $Z$ contains orbits of $T$ of every finite length,
and the restriction of $T$ on $Z$, $T|_Z:Z \to X$, is a compact
operator. 
\end{Thm}

The proof of Theorem \ref{main} is based on Theorem \ref{Thm3}. 
We first need to define the basis constant of a finite set of normalized
vectors of a Banach space in an analogous way of the definition of the 
basis constant of an infinite sequence.

\begin{Def}
Let $X$ be a Banach space, $n\in \N$, and $x_1,x_2,\ldots, x_n$ be
normalized elements of $X$.  We define the basis constant of 
$x_1,\ldots, x_n$  to be 
\[
\text{\em{bc}} \{x_1,\ldots, x_n\} := \sup\left\{|\alpha_1|,\ldots,
  |\alpha_n|\colon \ \left\|\sum^n_{i=1} \alpha_ix_i\right\| =
  1\right\}. 
\]
\end{Def}
Notice that
\[
\text{bc} \{x_1,\ldots, x_n\}^{-1} = \inf\left\{\left\| \sum^n_{i=1}
    \beta_ix_i\right\|\colon \ \max_{1\le i\le n} |\beta_i| =
  1\right\} ,
\]
and that $\text{bc} \{x_1,\ldots, x_n\} < \infty$ if and only if
$x_1,\ldots, x_n$ are linearly independent. 

Before stating Theorem~\ref{Thm3} recall that if $T$ is a quasi-nilpotent 
operator on a Banach space $X$, then for every $x \in X$ and $\eta >0$
there exists an increasing sequence $(i_n)_{n=1}^\infty$ in $\N$ such that 
$\| T^{i_n}x \| \le \eta \| T^{i_n -1}x \|$. Theorem~\ref{Thm3} asserts
that if $T$ is a strictly singular 1-1 operator on a Banach space $X$ then 
for arbitrarily small $\eta >0$ and $k \in \N$ there exists $x \in X$, 
$\| x \| =1$, such
that $\| T^i x \| \le \eta \| T^{i-1}x \|$ for $i=1,2, \ldots , k+1$, and 
moreover, the basis constant of $x, Tx/\| Tx \|, \ldots , T^k x/ \|T^k x \|$
does not exceed $1/ \sqrt{\eta}$. 

\begin{Thm}\label{Thm3}
Let $T$ be  a strictly singular 1-1 operator on a Banach space
$X$. Let $Y$ be an infinite dimensional subspace of $X$, $F$ be a
finite  codimensional subspace of $X$ and $k\in\N$.
Then there exists $\eta_0 \in (0,1)$ such that for every
$0<\eta\le \eta_0$ there exists $x\in Y$, $\|x\| = 1$ satisfying 
\begin{itemize}
\item[(a)] $T^{i-1}x\in F$ and $\|T^ix\|\le \eta \|T^{i-1}x\|$ for 
  $i=1,2,\ldots, k+1$, and
\item[(b)] $\text{\em{bc}} \left\{x, \frac{Tx}{\|Tx\|},\ldots,
    \frac{T^kx}{\|T^kx\|}\right\} \le \frac1{\sqrt{\eta}}$, 
\end{itemize}
(where $T^0$ denotes the identity operator on $X$).
\end{Thm}

We postpone the proof of Theorem \ref{Thm3}. 

\begin{proof}[Proof of Theorem \ref{main}]
Let $T$ be a strictly singular 1-1 operator on an infinite
dimensional Banach space $X$, and $Y$ be an infinite dimensional
subspace of $X$. Inductively for $n \in \N$ we construct a normalized sequence
$(z_n)_n \subset Y$, an increasing sequence of finite families $(z_j^*)_{j \in J_n}$ of
normalized functionals on $X$ (i.e. $(J_n)_n$ is an increasing sequence
of finite index sets), and a sequence $(\eta_n)_n \subset (0,1)$, as
follows: 

For $n=1$ apply Theorem \ref{Thm3} for
$F=X$ (set $J_1 = \emptyset$), $k=1$, to obtain
$\eta_1 < 1/2^6$ and $z_1 \in Y$, $\| z_1 \| =1$ such that 
\begin{equation} \label{eq62}
\| T^i z_1 \| < \eta_1 \| T^{i-1} z_1 \| \text{ for } i=1,2, 
\end{equation} 
and
\begin{equation} \label{eq63}
\text{bc} \{ z_1, \frac{Tz_1}{\| Tz_1 \| } \} < \frac{1}{\sqrt{
    \eta_1}}.
\end{equation} 
For the inductive step, assume that for $n \ge 2$, $(z_i)_{i=1}^{n-1} \subset Y$, 
$(z_j^*)_{j \in J_i}$ $(i=1, \ldots , n-1)$, and
$(\eta_i)_{i=1}^{n-1}$ have been constructed. Let $J_n$ be a finite
index set with $J_{n-1} \subseteq J_n$ 
and $(x_j^*)_{j \in J_n}$ be a set of normalized functionals on
$X$ such that 
\begin{equation} \label{eq64}
\begin{array}{c}
\text{for every } x \in \text{span} \{ T^iz_j: 1 \le j \le n-1 , 0 \le
i \le j \} \\
\text{there exists }j_0 \in J_n \text{ such that } | x_{j_0}^*(x)| \ge
\| x \| / 2.
\end{array}
\end{equation}
Apply Theorem \ref{Thm3} for 
$F= \cap_{j \in J_n} \text{ker} (x_j^*)$, and 
$k=n$, to
obtain $\eta_n < 1/(n^2 2^{2n+4})$ and $z_n \in Y$, $\| z_n \|=1$ such that 
\begin{equation} \label{eq65}
T^{i-1} z_n \in F \text{ and } \| T^i z_n \| < \eta_n \| T^{i-1} z_n \|
\text{ for } i=1,2, \ldots , n+1,
\end{equation}
and
\begin{equation} \label{eq66}
\text{bc} \{ z_n, \frac{Tz_n}{\| Tz_n \| }, \ldots , \frac{T^nz_n}{\|
  T^nz_n \| } \} < \frac{1}{\sqrt{\eta_n}}. 
\end{equation}
This finishes the induction.

Let $\widetilde{Z}= \text{span} \{ T^i z_n: n \in \N, 0 \le i \le n
\}$, and for $n \in \N$, let 
$Z_n= \text{span} \{ T^i z_n : 0 \le i \le n \}$.   Let $x \in
\widetilde{Z}$  with $\| x \| =1$ and write  $x= \sum_{n=1}^\infty x_n$ where 
$x_n \in Z_n$ for all $n \in \N$. We claim that 
\begin{equation} \label{eq67}
\| Tx_n \| < \frac{1}{2^n} \text{ for all } n \in \N.
\end{equation}
Indeed, write
$$
x=\sum_{n=1}^\infty \sum_{i=0}^n a_{i,n} \frac{T^i z_n}{\| T^i z_n \|}
\text{ and } x_n = \sum_{i=0}^n a_{i,n} \frac{T^i z_n}{\| T^i z_n \| }
\text{ for } n \in \N.
$$
Fix $n \in \N$ and set $\widetilde{x}_n= x_1+ x_2 + \cdots + x_n$. Let
$j_0 \in J_{n+1}$ such that 
\begin{align*}
\| \widetilde{x}_n \| & \le 2 | x^*_{j_0} (\widetilde{x}_n) |
\text{ (by (\ref{eq64}) for $n-1$ replaced by $n$)}\\
&= 2 | x^*_{j_0} (x) | \text{ (since for $n+1 \le m$, $J_{n+1} \subseteq J_m$ thus
by (\ref{eq65}), $x_m \in \text{ker} (x_{j_0}^*)$)}\\
& \le 2 \| x^*_{j_0} \| \| x \| =2.
\end{align*}
Thus $\| x_n \| = \| \widetilde{x}_n - \widetilde{x}_{n-1} \| \le 
\| \widetilde{x}_n \| + \| \widetilde{x}_{n-1} \| \le 4$ (where 
$\widetilde{x}_0=0$). Hence, by (\ref{eq63}) and (\ref{eq66}) we
obtain that 
\begin{equation} \label{eq68}
| a_{i,n} | \le 4 \text{bc} \{ \frac{T^iz_n}{ \| T^i z_n \| }: i=0,
\ldots , n \} \le \frac{4}{\sqrt{\eta_n}} \text{ for }i=0, \ldots , n.
\end{equation}
Therefore
\begin{align*}
\| Tx_n \| & = \| \sum_{i=0}^n a_{i,n} \frac{T^{i+1}z_n}{\| T^i z_n
  \|} \|  \le \sum_{i=0}^n | a_{i,n} | \frac{ \| T^{i+1}z_n \| }{\|
  T^iz_n\| }\\
&\le \sum_{i=0}^n \frac{4}{\sqrt{\eta_n}} \eta_n \text{ (by
  (\ref{eq62}), (\ref{eq65}), and (\ref{eq68}))}\\
&= 4 n \sqrt{\eta_n} < \frac{1}{2^n} \text{ (by the choice of
  $\eta_n$),}
\end{align*}
which finishes the proof of (\ref{eq67}). Let $Z$ to be the closure of
$\widetilde{Z}$. We claim that $T|_Z:Z \to X$ is a compact operator,
which will finish the proof of Theorem \ref{main}. Indeed, let
$(y_m)_m \subset \widetilde{Z}$ where  for all $m \in \N$ we have 
$\| y_m \| = 1$, and write $y_m= \sum_{n=1}^\infty y_{m,n}$ where
$y_{m,n} \in Z_n$ for all $n \in \N$. It suffices to prove that
$(Ty_m)_m$ has a Cauchy subsequence. Indeed, since $Z_n$ is finite
dimensional for all $n \in \N$, there exists $(y^1_m)_m$ a subsequence
of $(y_m)_m$ such that $(Ty^1_{m,1})_m$ is Cauchy. Let $(y^2_m)_m$ be
a subsequence of $(y^1_m)_m$ such that $(Ty^2_{m,2})_m$ is
Cauchy. Continue similarly, and let $\widetilde{y}_m = y^m_m$ and 
$\widetilde{y}_{m,n}= y^m_{m,n} $ for all $m,n \in \N$. Then for 
$m \in \N$ we have 
$\widetilde{y}_m= \sum_{n=1}^\infty \widetilde{y}_{m,n}$ where
$\widetilde{y}_{m,n} \in Z_n$ for all $n \in \N$. Also, for all 
$n,m \in \N$ with $n \le m$, $(\widetilde{y}_t)_{t \ge m}$
and $(\widetilde{y}_{t,n})_{t \ge m}$ are subsequences of 
$(y^m_t)_t$ and $(y^m_{t,n})_t$ respectively. Thus for all $n \in \N$,
$(T\widetilde{y}_{t,n})_{t \in \N}$ is a Cauchy sequence. We claim
that $(T\widetilde{y}_m)_m$ is a Cauchy sequence. Indeed, for
$\varepsilon >0$ let $m_0 \in \N$ such that $1/2^{m_0-1}< \varepsilon$
and let $m_1 \in \N$ such that 
\begin{equation} \label{eq69}
\| T \widetilde{y}_{s,n} - T \widetilde{y}_{t,n} \| <
\frac{\varepsilon}{2m_0} \text{ for all } s,t \ge m_1 \text{ and }
n=1,2, \ldots m_0.
\end{equation}
Thus for $s,t \ge m_1$ we have
\begin{align*}
\| T \widetilde{y}_s - T\widetilde{y}_t \| & = 
\| \sum_{n=1}^\infty T \widetilde{y}_{s,n} - T \widetilde{y}_{t,n} \| \\
& \leq \sum_{n=1}^{m_0} \| T \widetilde{y}_{s,n} - T \widetilde{y}_{t,n} \| 
+ \sum_{n=m_0 +1}^\infty \| T \widetilde{y}_{s,n} \| + 
\sum_{n=m_0+1}^\infty \| T \widetilde{y}_{t,n} \| \\
& < m_0 \frac{\varepsilon}{2m_0} + 2 \sum_{n=m_0 +1}^\infty 
\frac{1}{2^n} \text{ (by (\ref{eq67}) and (\ref{eq69}))}\\
&= \frac{\varepsilon}{2} + \frac{2}{2^{m_0}} < \varepsilon 
\text{ ( by the choice of $m_0$),}
\end{align*}
which proves that $(T\widetilde{y}_m)_m$ is a Cauchy sequence and
finishes the proof of Theorem~\ref{main}.
\end{proof}

For the proof of Theorem~\ref{Thm3} we need the next two results.

\begin{Lem}\label{lem1}
Let $T$ be a strictly singular 1-1 operator on an infinite
dimensional Banach  space $X$. Let $k\in\N$ and $\eta>0$. Then for
every  infinite dimensional subspace $Y$ of $X$ there exists an
infinite  dimensional subspace $Z$ of $Y$ such that for all $z\in Z$
and  for all $i=1,\ldots, k$ we have that
\[
\|T^iz\|\le \eta\|T^{i-1}z\|
\]
(where $T^0$ denotes the identity operator on $X$).
\end{Lem}

\begin{proof}
Let $T$ be a strictly singular 1-1 operator on an infinite dimensional
Banach space $X$, $k \in \N$ and $\eta >0$. 
We first prove the following \medskip

\noindent {\bf Claim:}\ For every infinite dimensional linear
submanifold (not necessarily closed) $W$
of $X$ there exists an infinite dimensional linear submanifold $Z$ of
$W$ such that $\|Tz\| \le \eta\|z\|$ for all $z\in Z$. \medskip

Indeed, since $W$ is infinite dimensional there exists a normalized
basic sequence $(z_i)_{i\in\N}$ in $W$ having basis constant at most
equal to 2, such that $\|Tz_i\|\le \eta/2^{i+2}$ for all $i\in\N$. Let
$Z = \text{span}\{z_i\colon i\in\N \}$ be the linear span of the
$z_i$'s. Then $Z$ is an infinite dimensional linear submanifold of
$W$. We now show that $Z$ satisfies the conclusion of the Claim. Let
$z\in Z$ and write $z$ in the form $z = \Sigma\lambda_iz_i$ for some
scalars $(\lambda_i)$ such that at most finitely many $\lambda_i$'s
are non-zero. Since the basis constant of $(z_i)_i$ is at most equal
to 2, we have that $|\lambda_i|\le 4\|z\|$ for all $i$.  Thus
\[
\|Tz\| = \left\|\sum_i \lambda_iTz_i\right\| \le \sum_i |\lambda_i|
\|Tz_i\| \le \sum_i 4\|z\| \frac\eta{2^{i+2}} =  \eta\|z\|
\]
which finishes the proof of the Claim.

Let $Y$ be an infinite dimensional subspace of $X$. Inductively for
$i=0,1, \ldots , k$, we define $Z_i$, a linear submanifold of
$X$, such that 
\begin{itemize}
\item[(a)] $Z_0$ is an infinite dimensional linear submanifold of $Y$
  and $Z_i$ is an infinite dimensional linear
  submanifold of  $T(Z_{i-1})$  for $i\ge 1$.
\item[(b)] $\|Tz\| \le \eta\|z\|$ for all $z\in Z_i$ and for all 
$i  \ge 0$.
\end{itemize}
 Indeed, since $Y$ is infinite dimensional, we obtain $Z_0$ by applying the
 above Claim for $W=Y$. Obviously (a) and (b) are satisfied for
 $i=0$. Assume that for some $i_0\in \{0,1,\ldots, k-1\}$, a linear
 submanifold $Z_{i_0}$ of $X$ has been constructed satisfying (a) and
 (b) for $i=i_0$. Since $T$ is 1-1 and $Z_{i_0}$ is infinite
 dimensional we have that $T(Z_{i_0})$ is an infinite dimensional
 linear submanifold of $X$ and we obtain $Z_{i_0+1}$ by applying the
 above Claim for $W = T(Z_{i_0})$. Obviously (a) and (b) are satisfied
 for $i=i_0+1$. This finishes the inductive construction of the
 $Z_i$'s. By (a) we obtain that $Z_k$ is an infinite dimensional
 linear submanifold of $T^k(Y)$. Let $W = T^{-k}(Z_k)$. Then $W$ is an
 infinite dimensional linear submanifold of $X$. Since $Z_k \subseteq
 T^k(Y)$ and $T$ is 1-1, we have that $W\subseteq Y$. By (a) we obtain
 that for $i=0,1,\ldots, k$ we have $Z_k \subseteq T^{k-i}Z_i$, hence 
$$
T^i W = T^i T^{-k} Z_k = T^{-(k-i)}Z_k \subseteq T^{-(k-i)}T^{k-i}Z_i
= Z_i
$$
(since $T$ is 1-1).  Thus by (b) we
 obtain that $\|T^iz\|\le \eta\|T^{i-1}z\|$ for all $z\in W$ and
 $i=1,2,\ldots, k$. Obviously, if $Z$ is the closure of $W$ then $Z$
 satisfies the statement of the  lemma.
\end{proof}

\begin{Cor}\label{cor2}
Let $T$ be a strictly singular 1-1 operator on an infinite
dimensional Banach space $X$. Let $k\in \N$, $\eta>0$ and $F$ be a
finite codimensional subspace of $X$. Then for every infinite
dimensional subspace $Y$ of $X$ there exists an infinite dimensional
subspace $Z$ of $Y$ such that for all $z\in Z$ and for all  $i=1,\ldots, k+1$
\[
T^{i-1}z\in F\quad \text{and}\quad \|T^iz\| \le \eta\|T^{i-1}z\|
\]
(where $T^0$ denotes the identity operator on $X$).
\end{Cor}

\begin{proof}
For any linear submanifold $W$ of $X$ and for any finite codimensional
subspace $F$ of $X$ we have that
\begin{equation}\label{eq1}
\dim(W/(F\cap W)) \le \dim(X/F) < \infty.
\end{equation}
Indeed for any $n>\dim(X/F)$ and for any $x_1,\ldots, x_n$ linear
independent vectors in $W\backslash(F\cap W)$ we have that there exist
scalars $\lambda_1,\ldots,\lambda_n$ with $(\lambda_1,\ldots,
\lambda_n) \ne (0,\ldots, 0)$ and $\sum\limits^n_{i=1} \lambda_i
x_i\in F$ (since $n>\dim(X/F)$). Thus $\sum\limits^n_{i=1}
\lambda_ix_i \in F\cap W$ which implies  (\ref{eq1}).

Let $R(T)$ denote the range of $T$. Apply
(\ref{eq1}) for $W=R(T)$ to obtain 
\begin{equation}\label{eq2}
\dim(R(T)/(R(T)\cap F)\le \dim(X/F) < \infty.
\end{equation}
Since $T$ is 1-1 we have that
\begin{equation}\label{eq3}
\dim(X/T^{-1}(F))\le \dim(R(T)/(R(T)\cap F)).
\end{equation}
Indeed, for any $n> \dim(R(T)/(R(T)\cap F))$ and for any 
$x_1,\ldots, x_n$ linear independent vectors of $X\backslash T^{-1}(F)$, we have
that $Tx_1,\ldots, Tx_n$ are linear independent vectors of
$R(T)\backslash T(T^{-1}(F)) = R(T)\backslash F$ (since $T$ is
1-1). Thus $Tx_1,\ldots, Tx_n\in R(T)\backslash (R(T)\cap F)$ and
since $n>\dim(R(T)/(R(T)\cap F))$, there exist scalars
$\lambda_1,\ldots,\lambda_n$ with 
$(\lambda_1,\ldots, \lambda_n)\ne (0,\ldots, 0)$ such that 
$\sum\limits^n_{i=1} \lambda_i Tx_i\in R(T)\cap F$. Therefore 
 $T\left(\sum\limits^n_{i=1} \lambda_ix_i\right) \in F$, and hence
 $\sum\limits^n_{i=1} \lambda_ix_i\in T^{-1}(F)$,
which proves (\ref{eq3}). By combining (\ref{eq2}) and (\ref{eq3}) we
obtain
\begin{equation}\label{eq4}
\dim(X/T^{-1}(F)) < \infty.
\end{equation}
By (\ref{eq4}) we have that 
\begin{equation}\label{eq5}
\dim(X/T^{-i}(F)) < \infty, \text{ for } i=1,2,\ldots, k.
\end{equation}
Thus $\dim(X/W_1) < \infty$ where $W_1 = F\cap T^{-1}(F)
\cap\cdots\cap T^{-k}(F)$. Therefore if we apply (\ref{eq1}) for $W=Y$
and $F=W_1$ we obtain 
\begin{equation}\label{eq6}
\dim(Y/Y\cap W_1)\le \dim(X/W_1) < \infty,
\end{equation}
and therefore $Y\cap W_1$ is infinite dimensional. 

Now use Lemma \ref{lem1}, replacing $Y$ by $Y\cap W_1$, to obtain an
infinite dimensional subspace $Z$ of $Y\cap W_1$ such that 
\[
\|T^iz\| \le \eta\|T^{i-1}z\|
\]
for all $z\in Z$ and $i=1,\ldots, k+1$.
Notice that for $z\in Z$ and $i=1,\ldots, k$ we have that $z\in W_1$
thus $T^{i-1}z\in F$. 
\end{proof}

Now we are ready to give the 

\begin{proof}[Proof of Theorem \ref{Thm3}]
We prove by induction on $k$ that for every infinite dimensional
subspace $Y$ of $X$, finite codimensional
subspace $F$ of $X$, $k\in \N$, function $f\colon (0,1)\to (0,1)$ such that
$f(\eta)\searrow 0$ as $\eta \searrow 0$, and for $i_0\in \{0\}\cup
\N$, there exists $\eta_0>0$ such that for every $0<\eta\le \eta_0$
there exists $x\in Y$, $\|x\| = 1$ satisfying 
\begin{itemize}
\item[(a$'$)] $T^{i-1}x\in F$ and $\|T^ix\| \le \eta\|T^{i-1}x\|$ for
  $i=1,2,\ldots, i_0+k+1$. 
\item[(b$'$)] $\text{bc} \left\{\frac{T^{i_0}x}{\|T^{i_0}x\|},
    \frac{T^{i_0+1}x}{\|T^{i_0+1}x\|},\ldots,
    \frac{T^{i_0+k}x}{\|T^{i_0+k}x\|}\right\} \le \frac1{f(\eta)}$. 
\end{itemize}
For $k=1$ let $Y, F,f$, and $i_0$ as above, and let 
$\eta_0 \in (0,1) $ satisfying
\begin{equation}\label{eq7}
f(\eta_0) < \frac1{62}.
\end{equation}
Let $0 < \eta\le \eta_0$. Apply Corollary \ref{cor2} for $k$ and
$\eta$ replaced by 
$i_0+1$ and $\eta /4$ respectively, to obtain an infinite dimensional
subspace $Z_1$ of $Y$  such
that for all $z\in Z_1$ and for $i=1,2,\ldots, i_0+2$ 
\begin{equation}\label{eq8}
T^{i-1}z\in F\quad \text{and}\quad \|T^iz\| \le \frac\eta4 \|T^{i-1}z\|.
\end{equation}
Let $x_1\in Z_1$ with $\|x_1\| = 1$. If $\text{bc} \{T^{i_0}
x_1/\|T^{i_0}x_1\|$, $T^{i_0+1}x_1/\|T^{i_0+1}x_1\|\} \le 1/f(\eta)$
then $x_1$ satisfies (a$'$) and (b$'$) for $k=1$, thus we may assume
that 
\begin{equation}\label{eq9}
\text{bc}\left\{\frac{T^{i_0}x_1}{\|T^{i_0}x_1\|},
  \frac{T^{i_0+1}x_1}{\|T^{i_0+1}x_1\|}\right\} > \frac1{f(\eta)}. 
\end{equation}
Let
\begin{equation}\label{eq10}
0 < \eta_2 \le \frac\eta4 \wedge \min_{1\le i\le i_0}
\frac{\|T^{i_0}x_1\|}{2 \|T^ix_1\|} \wedge 
\min_{i_0 < i  \le i_0 +2} \frac{\| T^i x_1 \| }{2 \| T^{i_0}x_1 \| }f(\eta). 
\end{equation}
Let $z^*_1,z^*_2\in X^*$, $\|z^*_1\| = \|z^*_2\| = 1$,
$z^*_1(T^{i_0}x_1) = \|T^{i_0}x_1\|$ and $z^*_2(T^{i_0+1}x_1) =
\|T^{i_0+1}x_1\|$. 
Since $\ker z^*_1 \cap \ker z^*_2$ is finite codimensional and $T$ is
1-1, by (\ref{eq5}) we have that 
\begin{equation}\label{eq11}
\dim(X/T^{-i_0}(\ker z^*_1\cap \ker z^*_2)) < \infty.
\end{equation}
Apply Corollary \ref{cor2} for $F,k$ and $\eta$ replaced by 
$F\cap T^{-i_0}(\ker z^*_1\cap \ker z^*_2)$, $i_0+2$ and 
$\eta_2$ respectively, to obtain an infinite dimensional subspace
$Z_2$ of $Y$  such that for all $z\in Z_2$ and
for all $i=1,2,\ldots, i_0+2$ 
\begin{equation}\label{eq12}
T^{i-1}z\in F\cap T^{-i_0}(\ker z^*_1\cap \ker z^*_2)\quad
\text{and}\quad \|T^iz\| \le \eta_2\|T^{i-1}z\|. 
\end{equation}
Let $x^*_1\in X^*$ with $\|x^*_1\| = x^*_1(x_1) = 1$ and let $x_2\in
Z_2\cap \ker x^*_1$ with 
\begin{equation}\label{eq12a}
\|T^{i_0}x_1\| = \|T^{i_0}x_2\|
\end{equation}
and let $x = (x_1+x_2)/\|x_1+x_2\|$. We will show that $x$ satisfies
(a$'$) and (b$'$) for $k=1$. 

We first show that (a$'$) is satisfied for $k=1$. Since
$x_1,Tx_1,\ldots, T^{i_0+1}x_1\in F$ (by (\ref{eq8})) and 
$x_2, Tx_2,\ldots, T^{i_0+1}x_2\in F$ (by (\ref{eq12})) we have that
$x,Tx,\ldots, T^{i_0+1}x\in F$.  Before showing that the norm estimate
of (a$'$) is satisfied, we need some preliminary estimates:
(\ref{eq12b})-(\ref{eq14a}). 

If  $1\le i < i_0$ (assuming that $2\le i_0$) then

\begin{alignat}{2}
\|T^ix_1\| &= \frac12\|T^{i_0}x_1\| \left(\frac{\|T^{i_0}x_1\|}{2\|T^ix_1\|}\right)^{-1}
&\quad &\phantom{\quad}\nonumber\\
&\le \frac12 \|T^{i_0}x_1\|\eta^{-1}_2 &\quad &\text{(by (\ref{eq10}))}\nonumber\\
&= \frac12 \|T^{i_0}x_2\|\eta^{-1}_2&\quad &\text{(by (\ref{eq12a}))}\nonumber\\
&\le \frac12 \eta^{i_0-i}_2\|T^ix_2\| \eta^{-1}_2 &\quad &\text{(by
  applying (\ref{eq12}) for $z=x_2$, $i_0-i$ times)}\nonumber\\ 
\label{eq12b}
&\le \frac12 \|T^ix_2\| &\quad &\text{(since $\eta_2\le 1$ by
  (\ref{eq10})).}
\end{alignat}
Thus, by (\ref{eq12b}), for $1\le i < i_0$ (assuming that
$2\le i_0$) we have
\begin{equation}\label{eq12c}
\|T^ix\| \|x_1+x_2\| = \|T^ix_1 + T^ix_2\| \le \|T^ix_1\| + \|T^ix_2\|
\le \frac32 \|T^ix_2\|
\end{equation}
and 
\begin{equation}\label{eq12d}
\|T^ix\| \|x_1+x_2\| = \|T^ix_1  + T^ix_2\|\ge \|T^ix_2\| -
\|T^ix_1\|\ge \frac12 \|T^ix_2\| .
\end{equation}
Also notice that 
\begin{equation}\label{eq12e}
\|T^{i_0}x\| \|x_1+x_2\| = \|T^{i_0}x_1 + T^{i_0}x_2\| \le
\|T^{i_0}x_1\| + \|T^{i_0}x_2\| = 2\|T^{i_0}x_1\| \text{ (by
  (\ref{eq12a})),}
\end{equation}
and
\begin{equation}\label{eq12f}
\|T^{i_0}x\| \|x_1+x_2\| = \|T^{i_0}x_1 + T^{i_0}x_2\| \ge
z^*_1(T^{i_0}x_1 + T^{i_0}x_2) = z^*_1 (T^{i_0}x_1) = \|T^{i_0}x_1\|
\end{equation}
(by (\ref{eq12}) for $z=x_2$ and $i=1$). Also for $i_0 < i \le i_0+2$
we have that by  applying (\ref{eq12}) for $z=x_2$, 
$i-i_0$ times, we obtain 
\begin{align}
\|T^ix_2\| &\le \eta^{i-i_0}_2 \|T^{i_0}x_2\|\nonumber\\
&\le \eta_2\|T^{i_0}x_1\| \quad \text{(by $\eta_2<1$ and
  (\ref{eq12a}))}\nonumber\\ 
&= \eta_2 \frac{2 \|T^{i_0}x_1\|}{\|T^{i}x_1\|} \frac12
\|T^ix_1\|\nonumber\\ 
\label{eq12g}
&\le \frac12 f(\eta) \|T^ix_1\|\quad \text{(by (\ref{eq10}))}\\
\label{eq12h}
&\le \frac12 \|T^ix_1\|.
\end{align}
Thus for $i_0<i\le i_0+2$ we have
\begin{align}
\|T^ix\| \|x_1+x_2\| &= \|T^ix_1 + T^ix_2\|\nonumber\\
&\le \|T^ix_1\| + \|T^ix_2\|\nonumber\\
\label{eq13}
&\le \frac32\|T^ix_1\|\quad \text{(by (\ref{eq12h})).}
\end{align}
Also for $i_0 < i\le i_0+2$ we have
\begin{align}
\|T^ix\| \|x_1+x_2\| &= \|T^ix_1 + T^ix_2\|\nonumber\\
&\ge \|T^ix_1\| - \|T^ix_2\|\nonumber\\
\label{eq14}
&\ge \frac12\|T^ix_1\|\quad \text{(by (\ref{eq12h})).}
\end{align}
Later in the course of this proof we will also need that
\begin{alignat}{2}
\|T^{i_0+1}x\| \|x_1+x_2\| &= \|T^{i_0+1}x_1 +
T^{i_0+1}x_2\|&&\nonumber\\ 
&\ge \|T^{i_0+1}x_1\| - \|T^{i_0+1}x_2\|&&\nonumber\\
&\ge \frac2{f(\eta)} \|T^{i_0+1}x_2\| - \|T^{i_0+1}x_2\| &\quad
&\text{(by (\ref{eq12g}))}\nonumber\\ 
&= \frac{2-f(\eta)}{f(\eta)} \|T^{i_0+1}x_2\|
&&\nonumber\\
\label{eq14a}
&\ge \frac1{f(\eta)} \|T^{i_0+1}x_2\|&\quad &\text{(since
  $f(\eta)<1$).}
\end{alignat} 
Finally we will show that for $1\le i\le i_0+2$ we have that $\|T^ix\|
\le \eta\|T^{i-1}x\|$. Indeed if $i=1$ then 
\begin{alignat}{2}
\|T^ix\| &= \frac1{\|x_1+x_2\|} \|Tx_1 + Tx_2\|
&&\nonumber\\
&\le \frac1{\|x_1+x_2\|} (\|Tx_1\| + \|Tx_2\|)
&&\nonumber\\
&\le \frac1{\|x_1+x_2\|} \left(\frac\eta4 \|x_1\| + \eta_2
  \|x_2\|\right)& \quad &\text{(by (\ref{eq8}) ($z=x_1$), and
  (\ref{eq12}) ($z=x_2$))} \nonumber\\ 
&\le \frac1{\|x_1+x_2\|} \left(\frac\eta4 \|x_1\| + \eta_2(\|x_1+x_2\|
  + \|x_1\|)\right)&&\nonumber\\ 
&= \frac1{\|x_1+x_2\|} \left(\frac\eta4 + \eta_2\right) x^*_1(x_1) +
\eta_2& \quad & \text{(by the choice of $x_1^*$)} \nonumber\\ 
&= \frac1{\|x_1+x_2\|} \left(\frac\eta4+\eta_2\right) x^*_1(x_1+x_2)
+ \eta_2&&\text{(since $x_2\in \ker x^*_1$)}\nonumber\\ 
&\le \frac\eta4 + 2\eta_2 &&\text{(since $\|x^*_1\| = 1$)}\nonumber \\
&\le \eta &&\left(\text{since $\eta_2<\frac\eta4$ by
    (\ref{eq10})}\right). \label{eq14b}
\end{alignat}
If $1<i<i_0$ (assuming that $3 \le i_0$) we have that
\begin{alignat}{2}
\frac{\|T^ix\|}{\|T^{i-1}x\|} &\le \frac{\frac32
  \|T^ix_2\|}{\frac12\|T^{i-1}x_2\|} &\quad &\text{(by (\ref{eq12c})
  and (\ref{eq12d}))}\nonumber\\ 
&< 3\eta_2 &\quad &\text{(by (\ref{eq12}))}\nonumber\\
\label{eq14c}
&< \eta&\quad &\text{(by (\ref{eq10})).}
\end{alignat}
If $i=i_0>1$ then
\begin{alignat}{2}
\frac{\|T^ix\|}{\|T^{i-1}x\|} &\le
\frac{2\|T^{i_0}x_1\|}{\frac12\|T^{i_0-1}x_2\|} &\quad &\text{(by
  (\ref{eq12e}) and (\ref{eq12d}))}\nonumber\\ 
&= 4\frac{\|T^{i_0}x_2\|}{\|T^{i_0-1}x_2\|}&\quad &\text{(by
  (\ref{eq12a}))}\nonumber\\ 
&< 4\eta_2&\quad &\text{(by (\ref{eq12}) for $z=x_2$ and $i=1$)}\nonumber\\
\label{eq14d}
&< \eta &\quad &\text{(by (\ref{eq10})).}
\end{alignat}
If $i_0 < i \le i_0+2$ then
\begin{alignat}{2}
\frac{\|T^ix\|}{\|T^{i-1}x\|} &\le \frac{\frac32\|T^ix_1\|}{\frac12
  \|T^{i-1}x_1\|} &\quad &\text{(by (\ref{eq13}) and
  (\ref{eq14}))}\nonumber\\ 
\label{eq14e}
&< \eta &\quad &\text{(by (\ref{eq8}) for $z=x_1$).} 
\end{alignat}
Now (\ref{eq14b}), (\ref{eq14c}), (\ref{eq14d}) and
(\ref{eq14e}) yield
that  for $1\le i \le i_0+2$ we have $\|T^ix\| \le \eta\|T^{i-1}x\|$, 
thus $x$ satisfies (a$'$) for $k=1$. Before proving 
that $x$ satisfies (b$'$) for $k=1$ we need some preliminary
estimates: (\ref{eq15})-(\ref{eq19}). 
By (\ref{eq9}) there exist
scalars $a_0,a_1$ with $\max(|a_0, |a_1|) = 1$ and $\|w\| < f(\eta)$
where 
\begin{equation}\label{eq15}
w = a_0 \frac{T^{i_0}x_1}{\|T^{i_0}x_1\|} + a_1 \frac{T^{i_0+1}x_1}{\|T^{i_0+1}x_1\|}.
\end{equation}
Therefore
$$
\left| |a_0| - |a_1|\right| = \left| \left\|a_0
    \frac{T^{i_0}x_1}{\|T^{i_0}x_1\|}\right\| - \left\|a_1
    \frac{T^{i_0+1}x_1}{\|T^{i_0+1}x_1\|}\right\|\right| \le \|w\| <
f(\eta). 
$$
Thus $1-f(\eta) \le |a_0|, |a_1| \le 1$ and hence 
\begin{equation}\label{eq16}
\frac{|a_1|}{|a_0|} \le \frac1{|a_0|} \le \frac1{1-f(\eta)}.
\end{equation}
Also by (\ref{eq15}) we obtain that
\[
T^{i_0}x_1 = \frac{\|T^{i_0}x_1\|}{a_0}w - \|T^{i_0}x_1\|
\frac{a_1}{a_0} \frac{T^{i_0+1}x_1}{\|T^{i_0+1}x_1\|} 
\]
and thus 
\begin{equation} \label{eq17}
T^{i_0}x = \frac1{\|x_1+x_2\|} \left(\frac{\|T^{i_0}x_1\|}{a_0} w -
  \|T^{i_0}x_1\| \frac{a_1}{a_0} \frac{T^{i_0+1}x_1}{\|T^{i_0+1}x_1\|}
  + T^{i_0}x_2\right). 
\end{equation} 
Let
\begin{equation}\label{eq18}
\widetilde w = T^{i_0}x + \frac{\|T^{i_0}x_1\|}{\|x_1+x_2\|}
\frac{a_1}{a_0} \frac{T^{i_0+1}x_1}{\|T^{i_0+1}x_1\|} -
\frac{T^{i_0}x_2}{\|x_1+x_2\|}. 
\end{equation}
Notice that (\ref{eq17}) and (\ref{eq18}) imply that 
$\widetilde{w}=( \| T^{i_0} x_1 \| / ( \| x_1 + x_2 \| a_0 ))w $ and hence 
\begin{align} \label{eq19}
\| \widetilde{w} \| & = 
\frac{\|T^{i_0}x_1\|}{\|x_1+x_2\||a_0|} \|w\| 
 \le \frac{\|T^{i_0}x_1\|}{\|x_1+x_2\|} \frac{f(\eta)}{1-f(\eta)}
\quad \text{(using (\ref{eq16}) and $\|w\| < f(\eta)$)} \nonumber\\ 
&\le 2f(\eta) \frac{\|T^{i_0}x_1\|}{\|x_1+x_2\|}\quad
\left(\text{since $\frac1{1-f(\eta)}<2$ by (\ref{eq7})}\right)
\nonumber \\ 
&= 2f(\eta) \frac{z^*_1(T^{i_0}x_1)}{\|x_1+x_2\|}
\quad \text{(by the choice of $z_1^*$)} 
\nonumber \\
&= 2f(\eta) \frac{z^*_1(T^{i_0}x_1 + T^{i_0}x_2)}{\|x_1+x_2\|}\quad
\text{(by (\ref{eq12}) for $i=1$ and $z=x_2$)} \nonumber \\ 
&\le 2f(\eta) \frac{\|T^{i_0}(x_1+x_2)\|}{\|x_1+x_2\|}
\quad \text{(since $\| z_1^* \| =1$)} \nonumber \\
&= 2f(\eta) \|T^{i_0}x\|.
\end{align}
Now we are ready to estimate the $\text{bc}\{T^{i_0}x/\|T^{i_0}x\|,
T^{i_0+1}x/\|T^{i_0+1}x\|\}$. Let scalars $A_0,A_1$ such that 
\[
\left\|A_0 \frac{T^{i_0}x}{\|T^{i_0}x\|} + A_1
  \frac{T^{i_0+1}x}{\|T^{i_0+1}x\|} \right\| = 1. 
\]
We want to estimate the $\max(|A_0|, |A_1|)$. By (\ref{eq18}) we have
\begin{align}
1 &= \left\|\frac{A_0}{\|T^{i_0}x\|} \left(\widetilde w -
    \frac{\|T^{i_0}x_1\|}{\|x_1+x_2\|} \frac{a_1}{a_0}
    \frac{T^{i_0+1}x_1}{\|T^{i_0+1}x_1\|} +
    \frac{T^{i_0}x_2}{\|x_1+x_2\|}\right) + A_1
  \frac{T^{i_0+1}x}{\|T^{i_0+1}x\|}\right\|\nonumber\\ 
&= \left\| \frac{A_0 \|T^{i_0}x_2\|}{\|T^{i_0}x\| \|x_1+x_2\|}
  \frac{T^{i_0}x_2}{\|T^{i_0}x_2\|} + \left(
    \frac{-A_0 \|T^{i_0}x_1\|}{\|T^{i_0}x\|\|x_1+x_2\|}
    \frac{a_1}{a_0} + 
\frac{A_1 \|T^{i_0+1}x_1\|}{\|T^{i_0+1}x\|\|x_1+x_2\|}\right)
  \frac{T^{i_0+1}x_1}{\|T^{i_0+1}x_1\|}
\right.\nonumber\\ 
&\quad\left.
   + \frac{ A_0 }{\|T^{i_0}x\|}  \widetilde w + 
 \frac{ A_1 T^{i_0+1}x_2 }{\|T^{i_0+1}x\|\|x_1+x_2\|}  \right\| 
\nonumber\\
&\ge \left\| \frac{A_0 \|T^{i_0}x_2\|}{\|T^{i_0}x\|\|x_1+x_2\|}
  \frac{T^{i_0}x_2}{\|T^{i_0}x_2\|} + \left(
    \frac{-A_0 \|T^{i_0}x_1\|}{\|T^{i_0}x\|\|x_1+x_2\|}
    \frac{a_1}{a_0}   +
\frac{A_1 \|T^{i_0+1}x_1\|}{\|T^{i_0+1}x\|\|x_1+x_2\|}\right)
  \frac{T^{i_0+1}x_1}{\|T^{i_0+1}x_1\|}\right\|\nonumber\\ 
\label{eq20}
&\quad -|A_0| 2f(\eta) - |A_1|f(\eta) \quad \text{(by the triangle
  inequality, (\ref{eq19}) and
  (\ref{eq14a})).} 
\end{align}
By (\ref{eq12}) for $i=1$ we have that $T^{i_0}x_2\in \ker z^*_2$ and
since $z^*_2(T^{i_0+1}x_1) = \|T^{i_0+1}x_1\|$ it is easy to see that
$\text{bc}\{T^{i_0}x_2/\|T^{i_0}x_2\|$, $T^{i_0+1}x_1/\|T^{i_0+1}x_1\|\} \le
2$. Thus (\ref{eq20}) implies that 
\begin{equation}\label{eq21}
\left|-\frac{A_0\|T^{i_0}x_1\|}{\|T^{i_0}x\| \|x_1+x_2\|}
  \frac{a_1}{a_0} + \frac{A_1\|T^{i_0+1}x_1\|}{\|T^{i_0+1}x\|
    \|x_1+x_2\|}\right| \le 2 + 4f(\eta) |A_0| + 2f(\eta) |A_1| 
\end{equation}
and
\begin{equation}\label{eq22}
\frac{|A_0| \|T^{i_0}x_2\|}{\|T^{i_0}x\|\|x_1+x_2\|} \le 2 +4f(\eta)
|A_0| + 2f(\eta)|A_1| .
\end{equation}
Notice that (\ref{eq22}) implies that
\begin{equation}\label{eq23}
|A_0| \le 4 + 8f(\eta) |A_0| + 4f(\eta)|A_1| ,
\end{equation}
since
\[
\frac{\|T^{i_0}x\| \|x_1+x_2\|}{\|T^{i_0}x_2\|} = \frac{\|T^{i_0}x_1 +
  T^{i_0}x_2\|}{\|T^{i_0}x_2\|} \le \frac{\|T^{i_0}x_1\| +
  \|T^{i_0}x_2\|}{\|T^{i_0}x_2\|} = 2 
\]
by (\ref{eq12a}). Also by (\ref{eq21}) we obtain
\[
\frac{|A_1| \|T^{i_0+1}x_1\|}{\|T^{i_0+1}x\| \|x_1+x_2\|} -
\frac{|A_0| \|T^{i_0}x_1\|}{\|T^{i_0}x\| \|x_1+x_2\|}
\frac{|a_1|}{|a_0|} \le2 + 4f(\eta) |A_0| + 2f(\eta)|A_1|. 
\]
Thus 
\begin{equation}\label{eq23a}
|A_1| \frac23 - |A_0| \frac1{1-f(\eta)} \le2 + 4f(\eta) |A_0| +
2f(\eta)|A_1|
\end{equation}
by (\ref{eq13}) for $i=i_0+1$, (\ref{eq16}) and
\begin{align*}
\frac{\|T^{i_0}x_1\|}{\|T^{i_0}x\|\|x_1+x_2\|} &=
\frac{\|T^{i_0}x_1\|}{\|T^{i_0}x_1 +T^{i_0}x_2\|} \le
\frac{\|T^{i_0}x_1\|}{z^*_1(T^{i_0}x_1 + T^{i_0}x_2)} \quad
\text{(since $\|z^*_1\|=1$)}\\ 
&= \frac{\|T^{i_0}x_1\|}{z^*_1(T^{i_0}x_1)} \quad \text{(since $x_2\in
  T^{-i_0}(\ker z^*_1)$ by (\ref{eq12}) for $i=1$ and $z=x_2$)}\\ 
&= 1\quad \text{(by the choice of $z^*_1$).}
\end{align*}
Notice that (\ref{eq23a}) implies that
\begin{equation}\label{eq24}
|A_1| \le 6 + \frac{28}5 |A_0|
\end{equation}
since $f(\eta) < 1/6$ by (\ref{eq7}). By substituting  (\ref{eq24})
into (\ref{eq23}) we obtain 
\begin{align*}
|A_0| &\le 4+8f(\eta) |A_0| + 4f(\eta) \left(6 + \frac{28}5 |A_0|\right)\\
&= 4+24f(\eta) + \frac{112}5 f(\eta)|A_0|\\
&\le 5 + \frac12 |A_0|\quad \left(\text{since $f(\eta) < \frac5{224}$
    by (\ref{eq7})}\right). 
\end{align*}
Thus $|A_0|\le 10$. Hence (\ref{eq24}) gives that $|A_1| \le 62$. Therefore 
\[
\text{bc}\left\{\frac{T^{i_0}x}{\|T^{i_0}x\|},
  \frac{T^{i_0+1}x}{\|T^{i_0+1}x\|}\right\} \le 62 \le \frac1{f(\eta)}
\quad \text{(by (\ref{eq7})).} 
\]
We now proceed to the inductive step. Assuming the inductive statement
for some integer $k$, let a finite codimensional subspace $F$ of $X$,
$f\colon (0,1)\to (0,1)$ with $f(\eta)\searrow 0$ as $\eta\searrow
0$ and $i_0\in\N\cup \{0\}$. By the inductive statement for $i_0,f$ and
$\eta$ replaced by $i_0+1$, $f^{1/4}$ and $\eta/4$ respectively, there
exists $\eta_1$ s.t.\ for $0<\eta<\eta_1$ there exists $x_1\in X$,
$\|x_1\| = 1$ 
\begin{equation}\label{eq25}
T^{i-1}x_1\in F \text{ and } \|T^ix_1\| \le \frac\eta4 \|T^{i-1}x_1\|
\text{ for } i=1,2,\ldots, (i_0+1) + k + 1 
\end{equation}
and
\begin{equation}\label{eq26}
\text{bc}\left\{\frac{T^{i_0+1}x_1}{\|T^{i_0+1}x_1\|},
  \frac{T^{i_0+2}x_1}{\|T^{i_0+2}x_1\|},\ldots,
  \frac{T^{i_0+1+k}x_1}{\|T^{i_0+1+k}x_1\|}\right\} \le
\frac1{f(\eta)^{1/4}}. 
\end{equation}
Let $\eta_0$ satisfying
\begin{equation}\label{eq27}
\eta_0 < \eta_1,\quad f(\eta_0) < \frac1{288^2},\quad f(\eta_0) <
\left(\frac1{144(k+1)}\right)^2 ,
\end{equation}
let $0 < \eta < \eta_0$ and let $x_1\in X$, $\|x_1\| = 1$ satisfying
(\ref{eq25}) and (\ref{eq26}). If 
\[
\text{bc}\left\{\frac{T^{i_0}x_1}{\|T^{i_0}x_1\|},
  \frac{T^{i_0+1}x_1}{\|T^{i_0+1}x_1\|},\ldots,
  \frac{T^{i_0+k+1}x_1}{\|T^{i_0+k+1}x_1\|}\right\} \le
\frac1{f(\eta)} 
\]
then $x_1$ satisfies the inductive step for $k$ replaced by
$k+1$. Thus we may assume that 
\begin{equation}\label{eq28}
\text{bc}\left\{\frac{T^{i_0}x_1}{\|T^{i_0}x_1\|},
  \frac{T^{i_0+1}x_1}{\|T^{i_0+1}x_1\|},\ldots,
  \frac{T^{i_0+k+1}x_1}{\|T^{i_0+k+1}x_1\|}\right\} >
\frac1{f(\eta)}. 
\end{equation}
Let
\begin{equation}\label{eq29}
0 < \eta_2 < \frac{\eta}{4} \wedge \min_{1 \le i \le i_0}
\frac{\| T^{i_0} x_1 \| }{2 \| T^i x_1 \| } \wedge 
\min_{i_0 < i \le i_0 +k +1} \frac{ \| T^i x_1 \| }{2 \| T^{i_0} x_1  \|}
f(\eta).
\end{equation}
Let $J \subset \{ 2,3, \ldots \}$ be a finite index set and 
$z^*_1, (z^*_j)_{j\in J}$ be norm 1 functionals such that
\begin{equation}\label{eq30}
z^*_1(T^{i_0}x_1) = \|T^{i_0}x_1\| ,
\end{equation}
and 
\begin{equation}\label{eq30a}
\text{for every } z\in \text{ span}\{T^{i_0+1}x_1,\ldots,
 T^{i_0+k+1}x_1\} \text{ there exists }j_0\in J \text{ with }|z^*_{j_0}(z)|\ge \frac12\|z\|.
\end{equation}
Since $T$ is 1-1  we obtain by (\ref{eq5}) that
$\dim(X/(T^{-i_0}\bigcap\limits_{j\in \{1\}\cup J} \ker z^*_j)) <
\infty$. Apply Corollary~\ref{cor2} for $F, k, \eta$
replaced by $F\cap T^{-i_0} \left(\bigcap\limits_{j\in \{1\}\cup J}
  \ker z^*_i\right)$, $i_0+k+2$, $\eta_2$ respectively, to obtain an
infinite dimensional subspace $Z$ of $Y$ such that for all
$z\in Z$ and for all $i=1,2,\ldots, i_0+k+2$ 
\begin{equation}\label{eq31}
T^{i-1}z \in F\cap T^{-i_0}\left( \bigcap_{j\in \{1\}\cup J} \ker
  z^*_j\right)\quad\text{and}\quad \|T^iz\| \le \eta_2
\|T^{i-1}z\|. 
\end{equation} 
Let $x^*_1\in X^*$, $\|x^*_1\| = 1 = x^*_1(x_1)$ and let $x_2\in Z\cap
\ker x^*_1$ with 
\begin{equation}\label{eq32}
\|T^{i_0}x_1\| = \|T^{i_0}x_2\|
\end{equation}
and let $x = (x_1+x_2)/\|x_1+x_2\|$. We will show that $x$ satisfies
the inductive statement for $k$ replaced by $k+1$. 

We first show that
$x$ satisfies (a$'$) for $k$ replaced by $k+1$. The proof is identical to the
verification of (a$'$) for $k=1$. The formulas (\ref{eq12g}),
(\ref{eq12h}), (\ref{eq13}), (\ref{eq14}), and (\ref{eq14e}) are valid
for $i_0  < i \le i_0+k+2$, and (\ref{eq14a}) is valid if $i_0+1$ is
replaced by any $i \in \{ i_0+1,\ldots, i_0+k+1\}$, and 
{\em this will be assumed in the rest of the proof when we refer to
  these formulas.}  

We now prove that (b$'$) is satisfied for $k$ replaced by $k+1$. By
(\ref{eq28}) there exist scalars $a_0,a_1,\ldots, a_{k+1}$ with
$\max(|a_0|, |a_1|,\ldots, |a_{k+1}|) = 1$ and $\|w\| < f(\eta)$ where

\begin{equation}\label{eq33}
w = \sum^{k+1}_{i=0} a_i \frac{T^{i_0+i}x_1}{\|T^{i_0+i}x_1\|}.
\end{equation}
We claim that
\begin{equation}\label{eq34}
|a_0| \ge \frac{f(\eta)^{1/4}}2.
\end{equation}
Indeed, if $|a_0| < f(\eta)^{1/4}/2$ then $\max(|a_1|,\ldots, |a_{k+1}|) = 1$ and
\begin{align*}
\left\|\sum^{k+1}_{i=1} a_i
  \frac{T^{i_0+1}x_1}{\|T^{i_0+i}x_1\|}\right\| &= \left\|w - a_0
  \frac{T^{i_0}x_1}{\|T^{i_0}x_1\|}\right\|\\ 
&\le \|w\| + |a_0|\\
&< f(\eta) + \frac{f(\eta)^{1/4}}{2}\\
&< f(\eta)^{1/4}\quad \text{(since $f(\eta) < 1/4$ by (\ref{eq27}))}
\end{align*}
which contradicts (\ref{eq26}). Thus (\ref{eq34}) is proved. By (\ref{eq33}) we obtain
\[
T^{i_0}x_1 = \frac{\|T^{i_0}x_1\|}{a_0} w - \sum^{k+1}_{i=1}
\frac{a_i}{a_0} \|T^{i_0}x_1\| \frac{T^{i_0+i}x_1}{\|T^{i_0+i}x_1\|} 
\]
and thus
\begin{equation} \label{pre51}
T^{i_0}x = \frac1{\|x_1+x_2\|} \left(\frac{\|T^{i_0}x_1\|}{a_0} w -
  \sum^{k+1}_{i=1} \frac{a_i}{a_0} \|T^{i_0}x_1\|
  \frac{T^{i_0+i}x_1}{\|T^{i_0+i}x_1\|} + T^{i_0}x_2\right). 
\end{equation}
Let
\begin{equation}\label{eq35}
\widetilde w = T^{i_0}x + \sum^{k+1}_{i=1} \frac{a_i}{a_0}
\frac{\|T^{i_0}x_1\|}{\|x_1+x_2\|}
\frac{T^{i_0+i}x_1}{\|T^{i_0+i}x_1\|} -
\frac{T^{i_0}x_2}{\|x_1+x_2\|}. 
\end{equation}
Notice that (\ref{pre51}) and (\ref{eq35}) imply that 
$\widetilde w = ( \| T^{i_0} x_1 \| / (\| x_1 + x_2 \| a_0))w$
and hence
\begin{align}
\| \widetilde w \| & = 
\frac{\|T^{i_0}x_1\|}{\|x_1+x_2\| |a_0|} \|w\| 
 < \frac{\|T^{i_0}x_1\|}{\|x_1+x_2\|} 2f(\eta)^{3/4} \quad
\text{(by $\|w\| \le f(\eta)$ and (\ref{eq34}))} \nonumber\\ 
&= \frac{z^*_1(T^{i_0}x_1)}{\|x_1+x_2\|} 2f(\eta)^{3/4} \quad
\text{(by (\ref{eq30}))} \nonumber\\ 
&= \frac{z^*_1(T^{i_0}x_1+T^{i_0}x_2)}{\|x_1+x_2\|} 2f(\eta)^{3/4}
\text{(by (\ref{eq31}) for $i=1$ and $z=x_2$)}\nonumber\\
&\le \frac{\|T^{i_0}(x_1+x_2)\|}{\|x_1+x_2\|} 2f(\eta)^{3/4} \quad
\text{(since $\|z^*_1\|=1$)}\nonumber\\ 
\label{eq36}
&= \|T^{i_0}x\| 2f(\eta)^{3/4}.
\end{align}
Now we are ready to estimate the
$\text{bc}\{T^{i_0+i}x_1/\|T^{i_0+i}x_1\|\colon \ i=0,1,\ldots, k+1\}$. Let
scalars\break $A_0,A_1,\ldots, A_{k+1}$ such that 
\[
\left\|\sum^{k+1}_{i=0} A_i \frac{T^{i_0+i}x}{\|T^{i_0+i}x\|}\right\| = 1.
\]
We want to estimate the $\max(|A_0|, |A_1|,\ldots, |A_{k+1}|)$. By
(\ref{eq35}) we have 
\begin{align}
1 &= \left\|\frac{A_0}{\|T^{i_0}x\|} \left(\widetilde w -
    \sum^{k+1}_{i=1} \frac{a_i}{a_0}
    \frac{\|T^{i_0}x_1\|}{\|x_1+x_2\|}
    \frac{T^{i_0+i}x_1}{\|T^{i_0+i}x_1\|} +
    \frac{T^{i_0}x_2}{\|x_1+x_2\|}\right) + \sum^{k+1}_{i=1} A_i
  \frac{T^{i_0+i}x}{\|T^{i_0+i}x\|}\right\|\nonumber\\ 
&= \left\|\frac{A_0 \|T^{i_0}x_2\|}{\|T^{i_0}x\|\|x_1+x_2\|}
  \frac{T^{i_0}x_2}{\|T^{i_0}x_2\|} + \sum^{k+1}_{i=1} 
 \left( \frac{a_i}{a_0} \frac{-A_0 \|T^{i_0}x_1\|}{\|T^{i_0}x\|
      \|x_1+x_2\|} + 
\frac{A_i \|T^{i_0+i}x_1\|}{\|T^{i_0+i}x\|\|x_1+x_2\|}\right)
  \frac{T^{i_0+i}x_1}{\|T^{i_0+i}x_1\|}
\right.\nonumber\\ 
&\quad \left.  
  + \frac{ A_0 }{\|T^{i_0}x\|}  \widetilde w  + \sum^{k+1}_{i=1}
 A_i \frac{ T^{i_0+i}x_2 }{\|T^{i_0+i}x\|\|x_1+x_2\|} \right\|
\nonumber\\
& \ge \left\| \frac{A_0 \|T^{i_0}x_2\|}{\|T^{i_0}x\| \|x_1+x_2\|}
  \frac{T^{i_0}x_2}{\|T^{i_0}x_2\|} + \sum^{k+1}_{i=1} \left(
    \frac{a_i}{a_0} \frac{-A_0 \|T^{i_0}x_1\|}{\|T^{i_0}x\|
      \|x_1+x_2\|} + 
    \frac{A_i \|T^{i_0+i}x_1\|}{\|T^{i_0+i}x\|\|x_1+x_2\|}\right)
  \frac{T^{i_0+i}x_1}{\|T^{i_0+i}x_1\|}\right\|\nonumber\\ 
\label{eq36a}
&\quad - |A_0| 2f(\eta)^{3/4} - \sum^{k+1}_{i=1} |A_i| f(\eta)\quad
\text{(by (\ref{eq36}) and (\ref{eq14a}); see the paragraph above (\ref{eq33}))}. 
\end{align}
By (\ref{eq31}) for $i=1$ and $z=x_2$ we obtain that 
$T^{i_0} x_2 \in \bigcap\limits_{j\in J} \ker z^*_j$ and by
(\ref{eq30a})  and (\ref{eq26}) it is easy to see that 
\[
\text{bc}\left\{\frac{T^{i_0}x_2}{\|T^{i_0}x_2\|},
  \frac{T^{i_0+i}x_1}{\|T^{i_0+i}x_1\|}\colon \ i=1,\ldots,
  k+1\right\} \le \frac2{f(\eta)^{1/4}} \vee 3. 
\]
Since $f(\eta) < \left(\frac23\right)^4$ (by (\ref{eq27})), we have
that $3 \le 2/f(\eta)^{1/4}$, hence 
\[
\text{bc}\left\{\frac{T^{i_0}x_2}{\|T^{i_0}x_2\|},
  \frac{T^{i_0+i}x_1}{\|T^{i_0+i}x_1\|}\colon \ i=1,\ldots,
  k+1\right\} \le \frac2{f(\eta)^{1/4}}. 
\]
Thus (\ref{eq36a}) implies that
\begin{equation}\label{eq37}
|A_0| \frac{\|T^{i_0}x_2\|}{\|T^{i_0}x\|\|x_1+x_2\|} \le
\frac2{f(\eta)^{1/4}} \left(1+2f(\eta)^{3/4} |A_0| + \sum^{k+1}_{j=1}
  |A_j| f(\eta)\right) ,
\end{equation}
and for $i=1,\ldots, k+1$
\begin{equation} \label{eq38}
\left|\frac{a_i}{a_0}
  \frac{-A_0 \|T^{i_0}x_1\|}{\|T^{i_0}x\|\|x_1+x_2\|} + 
  \frac{A_i \|T^{i_0+i}x_1\|}{\|T^{i_0+i}x\|\|x_1+x_2\|}\right| 
 \le \frac2{f(\eta)^{1/4}} \left(1+2f(\eta)^{\frac34} |A_0| +
  \sum^{k+1}_{j=1} |A_j| f(\eta)\right). 
\end{equation}
Since
\[
\frac{\|T^{i_0}x\|\|x_1+x_2\|}{\|T^{i_0}x_2\|} = \frac{\|T^{i_0}x_1 +
  T^{i_0}x_2\|}{\|T^{i_0}x_2\|} \le \frac{||T^{i_0}x_1\| +
  \|T^{i_0}x_2\|}{\|T^{i_0}x_2\|} = 2 \quad \text{(by (\ref{eq32})),} 
\]
we have that (\ref{eq37}) implies 
\begin{equation}\label{eq39}
|A_0| \le \frac4{f(\eta)^{1/4}} + 8f(\eta)^{1/2} |A_0| +
4 \sum^{k+1}_{j=1} |A_j| f(\eta)^{3/4} .
\end{equation}
Notice also that (\ref{eq38}) implies that for $i=1,\ldots, k+1$
\[
|A_i| \frac{\|T^{i_0+i}x_1\|}{\|T^{i_0+i}x\| \|x_1+x_2\|} - |A_0|
\frac{|a_i|}{|a_0|} \frac{\|T^{i_0}x_1\|}{\|T^{i_0}x\| \|x_1+x_2\|}
\le \frac2{f(\eta)^{1/4}} + 4f(\eta)^{\frac12} |A_0| +
2 \sum^{k+1}_{j=1} |A_j| f(\eta)^{\frac34}. 
\]
Thus
\begin{equation}\label{eq40}
|A_i| \frac23 - |A_0| \frac2{f(\eta)^{1/4}} \le \frac2{f(\eta)^{1/4}}
+ 4f(\eta)^{\frac12} |A_0| + 2 \sum^{k+1}_{j=1} |A_j| f(\eta)^{\frac34} 
\end{equation}
by (\ref{eq13}) (see the paragraph above (\ref{eq33})), (\ref{eq34}) and
\begin{align*}
\frac{\|T^{i_0}x_1\|}{\|T^{i_0}x\| \|x_1+x_2\|} &=
\frac{\|T^{i_0}x_1\|}{\|T^{i_0}x_1 + T^{i_0}x_2\|} \le
\frac{\|T^{i_0}x_1\|}{| z^*_1(T^{i_0}x_1 + T^{i_0}x_2) |} \quad
\text{(since $\|z^*_1\| = 1$)}\\ 
&= \frac{\|T^{i_0}x_1\|}{ | z^*_1(T^{i_0}x_1) |}\quad \text{(since
  $T^{i_0}x_2\in \ker z^*_1$ by (\ref{eq31}) for $i=1$ and $z=x_2$)}\\ 
&= 1\quad \text{(by (\ref{eq30})).}
\end{align*}
For $i=1,\ldots, k+1$ rewrite (\ref{eq40}) as
\[
|A_i| \left(\frac23 - 2 f(\eta)^{3/4}\right) \le \frac2{f(\eta)^{1/4}} +
\left(4f(\eta)^{1/2} + \frac2{f(\eta)^{1/4}}\right) |A_0| +
\sum^{k+1}_{\stackrel{\scriptstyle j=1}{\scriptstyle j\ne i}} |A_j|
f(\eta)^{3/4}. 
\]
Thus, since $f(\eta) < \left(\frac16\right)^{4/3} \wedge 
  \left(\frac14\right)^{1/2}$ (by (\ref{eq27})), we obtain 
\[
|A_i| \frac13 \le \frac2{f(\eta)^{1/4}} + \left(1 +
  \frac2{f(\eta)^{1/4}}\right) |A_0| +
\sum^{k+1}_{\stackrel{\scriptstyle j=1}{\scriptstyle j\ne i}} |A_j|
f(\eta)^{3/4}. 
\]
Hence, since $1\le 1/f(\eta)^{1/4}$, we obtain that for $i=1,\ldots, k+1$
\begin{equation}\label{eq41}
|A_i| \le \frac6{f(\eta)^{1/4}} + \frac9{f(\eta)^{1/4}} |A_0| + 3
\sum^{k+1}_{\stackrel{\scriptstyle j=1}{\scriptstyle j\ne i}} |A_j|
f(\eta)^{3/4} .
\end{equation}
By substituting (\ref{eq39}) in
(\ref{eq41}) we obtain that for $i=1,\ldots, k+1$, 
\begin{equation}\label{eq42}
|A_i|\le \frac6{f(\eta)^{1/4}} + \frac{36}{f(\eta)^{1/2}} + 72 f(\eta)^{\frac{1}{4}}
|A_0| + 36 \sum^{k+1}_{j=1} |A_j| f(\eta)^{1/2} + 3
\sum^{k+1}_{\stackrel{\scriptstyle j=1}{\scriptstyle j\ne i}} |A_j|
f(\eta)^{3/4}. 
\end{equation}
We claim that (\ref{eq39}) and (\ref{eq42}) imply that $\max\{|A_i|\colon \ 0
\le i \le k+1\} \le 1/f(\eta)$ which finishes the proof. Indeed, if
$\max\{|A_i|\colon \ 0 \le i \le k+1\} = |A_0|$ then (\ref{eq39})
implies that 
\begin{align*}
|A_0| &\le \frac4{f(\eta)^{1/4}} + 8f(\eta)^{1/2} |A_0| + 4 (k+1) |A_0| f(\eta)^{3/4}\\
&\le \frac4{f(\eta)^{1/4}} + \frac13 |A_0| + \frac13 |A_0| \quad
\left(\text{since  $f(\eta) <
    \left(\frac1{24}\right)^2 \wedge 
    \left(\frac1{12(k+1)}\right)^{\frac43}$ by }(\ref{eq27})\right) 
\end{align*}
thus
\begin{equation}\label{eq43}
|A_0| \le \frac{12}{f(\eta)^{1/4}} < \frac1{f(\eta)} \quad
\left(\text{since  $f(\eta) <
    \left(\frac1{12}\right)^{4/3}$ by} (\ref{eq27}) \right). 
\end{equation}
Similarly, if there exists $\ell\in\{1,\ldots, k+1\}$ such that
$\max\{|A_i|\colon \ 0\le i \le k+1\} = |A_\ell|$ then (\ref{eq42})
for $i=\ell$ implies that 
\begin{align*}
|A_\ell| &\le \frac6{f(\eta)^{1/4}} + \frac{36}{f(\eta)^{1/2}} +
72f(\eta)^{\frac{1}{4}}|A_\ell| + 36(k+1) f(\eta)^{1/2} |A_\ell| + 3kf(\eta)^{3/4}
|A_\ell|\\ 
&\le \frac{42}{f(\eta)^{1/2}} + \frac14 |A_\ell| + \frac14 |A_\ell| +
\frac14 |A_\ell| 
\end{align*}
(since $1/f(\eta)^{1/4} \le 1/f(\eta)^{1/2}$ 
and  $f(\eta) < \frac1{288^4} \wedge \left(\frac1{144(k+1)}\right)^2 $ by
    (\ref{eq27})).  Hence
\begin{equation}\label{eq44}
|A_\ell| \le \frac{168}{f(\eta)^{1/2}} \le \frac1{f(\eta)} \quad
\left(\text{since $f(\eta) < \frac1{168^2}$ by (\ref{eq27})}\right). 
\end{equation}
By (\ref{eq43}) and (\ref{eq44}) we have that $\max\{|A_i|\colon \ 0\le i
\le k+1\} \le 1/f(\eta)$ which finishes the proof. 
\end{proof}

\vspace{.2in}
\scriptsize{
\noindent 
Department of Mathematics, University of South Carolina, Columbia, SC
29208. 
giorgis@math.sc.edu

\noindent 
Department of Mathematics, Kent State University, Kent,  OH 44240. 
enflo@mcs.kent.edu
}

\end{document}